\documentclass[a4paper,12pt]{article}
\usepackage{amsmath}
\usepackage{ulem}
\usepackage{calc}
\usepackage{vmargin}
\usepackage{amstext}
\usepackage{graphicx}
\usepackage{amsfonts}
\usepackage{amssymb}

\begin{document}

	\title{		{\bf Exponential map of germs of vector fields 
			 } }
	\author{ \bf O. V. Kaptsov
		\\ Federal Research Center for\\
		  Information and Computational Technologies
				\\ Novosibirsk, Russia
		\\ E-mail: kaptsov@mail.ru}
	
	\date{}
	\maketitle
	
\begin{center} {\bf Abstract.}  \end{center}

\noindent
In this paper we introduce an exponential map of the algebra of germs of vector fields into the group of germs of diffeomorphisms at zero.
It is shown that this mapping is not a bijection. A brief review of the key results of the analytic iteration problem is given.


The exponential map establishes a correspondence between finite-dimensional Lie algebras and Lie groups and is a local diffeomorphism \cite{Hall}. However, infinite-dimensional Lie algebras arise in applications. Therefore, it is necessary to study the connections between infinite-dimensional objects.

If we consider the group of diffeomorphisms of a manifold $M$ and the algebra of smooth vector fields $V(M)$ on $M$, then the exponential map
$$Exp: V(M) \rightarrow Diff(M)$$
associates a vector field with its flow at time $t=1$. This map is no longer a local bijection \cite{Milnor1}. To avoid global constructions, we will consider neighborhoods of zero in $\mathbb{R}^n$ and diffeomorphisms of these neighborhoods that leave zero fixed. Each such diffeomorphism defines a germ $u$ such that $u(0)=0$. The set of such germs forms a group $\text{diff}(\mathbb{R}^n, 0)$.

Similarly, consider smooth vector fields $v$ defined in neighborhoods of zero in $\mathbb{R}^n$, such that $v(0)=0$. Each vector field $v$ defines a germ at zero. The set of such germs of vector fields is denoted by $V(\mathbb{R}^n, 0)$, and it forms a Lie algebra.

We want to define the exponential map
$$\exp: V(\mathbb{R}^n, 0) \to \text{diff}(\mathbb{R}^n, 0)$$
and introduce one-parameter groups. Let $U$ be a neighborhood of zero in $\mathbb{R}^n$, $v$ be a smooth vector field on $U$ with $v(0)=0$. Then $v$ defines a system of ordinary differential equations
\begin{equation} \label{eq1}
	\frac{dy}{dt} = v(y), \quad y(0) = x.
\end{equation}
By the theorem on the existence and uniqueness of solutions to systems of ordinary differential equations \cite{Zwillinger}, for any $T \in \mathbb{R}$, there exists a ball $B$ centered at zero such that the solution of system \eqref{eq1} will exist on the interval $[-T, T]$ if the initial data lies in this ball. This is due to the fact that for any $\epsilon > 0$, we can choose a ball $B$ such that
$$\max_{y \in B} |v(y)| < \epsilon .$$
Therefore, the local flow $\varphi_t$ generated by the vector field $v$ in the ball $B$ is defined for any $t \in [-T, T]$. Thus, the germ of this vector field $v$ generates a germ $\varphi^t$ of the flow defined for all $t \in \mathbb{R}$, and a one-parameter group
$$\varphi^t: \mathbb{R} \to \text{diff}(\mathbb{R}^n, 0) .$$
Then, the exponential map
$$\exp: V(\mathbb{R}^n, 0) \to \text{diff}(\mathbb{R}^n, 0)$$
associates the germ of a vector field $v$ with the germ of the flow $\varphi^t$ at $t=1$.

The problem of constructing a germ of a flow $\varphi^t$ that becomes a given germ $u$ at $t=1$ can be formulated as the following conditions:
\begin{equation} \label{group}
	\varphi^t \circ \varphi^s = \varphi^{t+s}, \quad \varphi^t(0) = 0, \quad \varphi^0(x) = x, \quad \varphi^1 = u .
\end{equation}
We will additionally assume that the Jacobian matrix $J=(\frac{\partial u}{\partial x})$ satisfies the inequality
\begin{equation} \label{neq1}
	\left|J(0) - E\right| < 1,
\end{equation}
where $E$ is the identity matrix. This condition means that the germ $u$ is close to the identity germ $I(x)=x$, which is the identity element of the group $diff(\mathbb{R}^n, 0)$.

The problem \eqref{group}, in the one-dimensional case, is called the problem of analytic iteration. We will restrict ourselves to considering germs of analytic diffeomorphisms and discuss the conditions for the existence of a solution to the analytic iteration problem. The Jacobian matrix in the one-dimensional case is a number $\lambda=u^{\prime}(0)$, which is called the multiplier.

As shown by Koenigs \cite{Milnor2}, if $0<\lambda<1$, then the germ $u$ is conjugate to its linear part, i.e., there exists a germ $f \in \text{diff}(\mathbb{R}^n, 0)$ such that
$$f \circ u \circ f^{-1} = \lambda \cdot I .$$
Consequently, $u$ can be represented as
$$u = f^{-1}(\lambda f).$$
Then the desired flow $\varphi^t$ has the form
$$\varphi^t = f^{-1}(\lambda^t f).$$
If the multiplier $\lambda$ is equal to $1$, then an analytic flow satisfying conditions \eqref{group} may not exist, although a solution in the form of a formal power series always exists. In particular, an analytic solution does not exist for the germ $u=e^x-1$.

If we consider the field of complex numbers $\mathbb{C}$ instead of the field of real numbers $\mathbb{R}$, then the problem of local conjugacy will have a solution if the multiplier $\lambda$ satisfies the inequality $|\lambda| \neq 1$. This follows from Koenigs' conjugacy theorem. The case $|\lambda| = 1$ is much more complicated.

Example 1. Let
$$u(z) = e^{i\pi/m} z + z^{2m+1}, \quad m \geq 1 ,$$
where $m$ is a natural number. For this $u$, the problem of analytic iteration has no solution even in the class of formal power series. Indeed, if such a flow $\varphi^t$ exists, then for $t = \frac{1}{2}$, the equality
$$\varphi^{1/2} \circ \varphi^{1/2}  = u$$
must hold. It is easy to see that there is no power series of the form
\begin{equation} \label{g=}
	g(z) = \sum_{i=1}^\infty c_i z^i
\end{equation}
such that
\begin{equation} \label{g g}
	g(g(z)) = e^{i\pi/m} z + z^{2m+1} .
\end{equation}
Substituting the series \eqref{g=} into the left-hand side of \eqref{g g}, collecting like terms for $z, z^2, ..., z^{2m+1}$, and comparing with the right-hand side of \eqref{g g}. Collecting like terms for $z$, we obtain the equation for $c_1$:
\begin{equation} \label{c1=}
	c_1^2 = e^{i\pi/m} .
\end{equation}

Collecting like terms at $z^2, ..., z^{2m}$, we find $c_2 = c_3 = ... = c_{2m} = 0$. Finally, collecting like terms at $z^{2m+1}$, we have
$$ c_{2m+1}c_1 (c_1^{2m}+1) = a\neq 0 . $$
The last expression contradicts (\ref{c1=}). If $m>4$, then $| e^{i\pi/m} - 1| = 2-\sqrt{2}< 1 $. Hence, condition (\ref{neq1}) is satisfied.

{\it Remark.} In the previous example, the multiplier $\lambda = e^{i\pi/m}$ satisfies the resonance condition: $\lambda= \lambda^{2m+1} $. The presence of the resonant term $z^{2m+1}$ prevents the existence of a solution to problem (\ref{group}).

Let us return to the general case. Let a matrix $A$ of size $n \times n$ be given. The set of eigenvalues $\lambda_1, ..., \lambda_n$ of the matrix $A$ is called resonant if there exist numbers $s, m_1,\dots,m_n$ such that
$$ \lambda_s = \lambda_1^{m_1} ... \lambda_n^{m_n}, \quad m_i \geq 0, \quad \sum m_i \geq 2.$$

The following Poincaré theorem \cite{Arnold} holds.\\
{\bf Theorem.} If the eigenvalues of the Jacobi matrix $J(0)$ of a germ $u \in \text{diff}(\mathbb{R}^n, 0)$ are less than one in modulus and there are no resonances, then there exists a germ $f \in \text{diff}(\mathbb{R}^n, 0)$ such that $u$ is conjugate to the linear germ $I(x)=x$, i.e.,
$$f \circ u \circ f^{-1}(x) = I(x).$$

\noindent
{\bf Corollary.} If the conditions of the Poincaré theorem are satisfied, then the problem (\ref{group}) has the solution $$\varphi^t = f^{-1}(J^t(0) f), \quad J^t(0) = e^{t \ln J(0)}.$$

{\it Remark.} A generalization of the Poincaré theorem is Siegel's theorem \cite{Arnold}. In the case of resonances, the germ $u$ is reduced to a normal form by a formal transformation \cite{Arnold}.

Let us give an example of a germ preserving area for which problem (\ref{group}) is not solvable even in the class of formal power series. Consider the rotation matrix
$$M = \begin{pmatrix} \cos \alpha & -\sin \alpha \\ \sin \alpha & \cos \alpha \end{pmatrix}, \quad \alpha = \pi/m, \qquad m\geq 1$$
and the polynomial germ
$$\psi = x + y^{2m}, \quad y = y.$$
Obviously, the germ $\psi$ preserves area. Hence, the composition of the rotation and the germ $\psi$ also preserves area. This composition, using complex variables $z = x + iy$, $\bar{z} = x - iy$, is written as
$$u = e^{i\pi/m} \left(z + \left(\frac{z - \bar{z}}{2i}\right)^{2m+1}\right).$$
It is easy to see that there is no series
\begin{equation} \label{gz=}
	g(z, \bar{z}) = \sum_{i, j \geq 1} c_{ij} z^i \bar{z}^j
\end{equation}
such that
\begin{equation} \label{gg=u}
	g \circ g = u .
\end{equation}
Substituting (\ref{gz=}) into the left-hand side of (\ref{gg=u}), collecting like terms at $z^i \bar{z}^j$, and comparing with the right-hand side of (\ref{gg=u}), we successively obtain
$$c_{10} = e^{i\pi/m}, \quad c_{ij} = 0 \quad \text{for} \quad 1 < i + j < 2m+1.$$
Collecting terms with $z^{2m+1}$, we find
$$\exp(i\pi/m) (g_1)^{2m+1} = c_{2m+1, 0} + c_{1, 0} (1 + c_1) = 0.$$
This contradiction proves the absence of a formal series. This example is taken from the work \cite{Kaptsov}. It shows that the statement of Yu. Moser (\cite{Moser} pp. 24-34) on the existence of a formal solution to problem (2) is incorrect.

It is worth noting the important work \cite{Lewis} devoted to the formal solvability of problem (2). Historical studies devoted to the one-dimensional problem of analytic iteration are presented in \cite{Kaptsov}. However, many questions related to the problem of inverting the exponential mapping remain open at present.

Note that the approach using germs of vector fields and transformation groups may be useful in the theory of differential equations.

 \textbf{Acknowledgements.} 

The research was carried out within the state assignment of Ministry of Science and Higher Education of the Russian Federation for Federal Research Center for Information and Computational Technologies.

This work is supported by the Krasnoyarsk Mathematical Center and financed by the Ministry of Science and Higher Education of the Russian Federation in the framework of the establishment and development of regional Centers for Mathematics Research and Education (Agreement No. 075-02-2025-1606).

\end{document}